\documentclass[12pt,leqno]{article} 

\usepackage{amsmath} \usepackage{amscd} 
\usepackage{amssymb} 
\usepackage{amsfonts}
\usepackage{amsthm}
\usepackage{makeidx} 
\usepackage{verbatim}
\newlength{\fixboxwidth}
\newtheorem{thm}{Theorem}
\newtheorem{cor}{Corollary}
\newtheorem{rem}{Remark}
\newtheorem{definition}{Definition}
\newtheorem{lemma}{Lemma}

\def\wt{\widetilde }

\def\e{\varepsilon } 
\def\epsilon{\varepsilon }
 
\def\rho{\varrho }

\def\phi{\varphi }

\def\({\biggl( }
\def\){\biggr) }

\def\span{{\rm span\,}}

\def\sq2{{\sqrt{2}}}

\def\A{{\mathcal A}}
  
\def\lin{{\rm lin}}
\def\non{{\rm non}}
\def\cont{{\rm cont}}

\newcount\refnum
\def\refx{\smallskip \global\advance\refnum by 1 {[\the\refnum ] \ }}

\def\B{{\cal B}}
\def\L{{\cal L}}
\def\N{{\cal N}}

\def\C{{\cal C}}
\def\D{{\cal D}}

\def\Nb{{\mathbb N}}
\def\R{{\mathbb R}} 
\def\Rd{{\mathbb R}^d}

\def\supp{{\rm supp \, }}

\def\lsim{\raisebox{-1ex}{$~\stackrel{\textstyle <}{\sim}~$}}

\title{Optimal Approximation of Elliptic Problems by Linear
and Nonlinear Mappings I} 

\author{Stephan Dahlke\thanks{The work of this author 
has been supported through the
European Union's Human Potential Programme, under contract
HPRN-CT-2002-00285 (HASSIP), and through DFG, 
Grant  Da 360/4-2.}, Erich Novak, Winfried Sickel} 

\begin{document}

\maketitle

\begin{abstract}  
We study the optimal approximation of the solution of an operator 
equation $ \A(u) = f $ 
by linear mappings of rank $n$ and compare this with the best
$n$-term approximation with respect to an optimal Riesz basis. 
We consider 
worst case errors, where $f$ is an element of the unit ball of 
a Hilbert space.
We apply our results to boundary value problems for 
elliptic PDEs that are given by an isomorphism 
$ 
\A : H^s_0 (\Omega) \to H^{-s} (\Omega) , 
$ 
where $s >0$ and $\Omega$ is an arbitrary bounded Lipschitz 
domain in $\R^d$. 
We prove that approximation by linear mappings is as good 
as the best $n$-term approximation with respect to an optimal 
Riesz basis. 
We discuss why nonlinear approximation still is
important for the approximation of elliptic problems. 
\end{abstract}

%partial differential equation, boundary value problems, miscellaneous


\section{Introduction}

We study the optimal approximation of the solution of an operator 
equation
\begin{equation}     \label{eq01} 
\A(u) = f , 
\end{equation} 
where $\A$ is a linear operator 
\begin{equation}     \label{eq02}
\A: H  \to G
\end{equation} 
from a Hilbert space $H$ to another Hilbert space $G$.
We always assume that $\A$ is boundedly invertible, hence 
\eqref{eq01} has a unique solution for any $f \in G$.  
We have in mind, for example, the more specific situation of an 
elliptic operator equation,  which is given as follows. 
Assume that $\Omega \subset \R^d$ is a bounded Lipschitz
domain and assume that 
\begin{equation}     \label{eq04} 
\A : H^s_0 (\Omega) \to H^{-s} (\Omega) 
\end{equation} 
is an isomorphism, where $s > 0$. 
A standard case (for second order elliptic boundary value problems for
PDEs) is $s=1$, but also other values of $s$ are of
interest.
For this situation we take $H=H^s_0(\Omega)$ and $G = H^{-s} (\Omega)$. 
Since $\A$ is boundedly invertible, the inverse mapping 
$S: G \to H$ is well defined. This mapping is sometimes 
called the solution operator --- in particular if we want 
to compute the solution $u=S(f)$ from the given 
right-hand side $\A(u)=f$. 

Let $F$ be a specified normed (or quasi-normed) subspace of $G$. 
We use linear and  nonlinear mappings $S_n$   for
approximating the solution $u=\A^{-1}(f)$  for $f \in F$.
Let us consider the worst case error 
$$ 
e(S_n, F,H) = 
\sup_{\Vert f \Vert_F \le 1} \Vert \A^{-1}(f)- S_n(f) \Vert_H . 
$$ 
%  where $F$ is a normed (or quasi-normed) space, $F \subset G$. 
For a given basis
$\B = \{ h_i \mid i \in \Nb \}$ of $H$ 
we consider the class $\N_n (\B)$ of all (linear or 
nonlinear) mappings of the form 
$$  
S_n(f) = \sum_{k=1}^n c_k \,  h_{i_k} ,      
$$ 
where the $c_k$ and the $i_k$ depend in an arbitrary way on $f$.
We also allow the basis $\B$ to be chosen in a nearly arbitrary way. 
Then the nonlinear widths
$e_{n,C}^\non (S,F,H)$ are given by
$$ 
e_{n,C}^\non (S, F , H) = \inf_{\B \in \B_C} \inf_{S_n \in \N_n(\B)} 
e(S_n, F,H). 
$$ 
Here $\B_C$ denotes a set of Riesz bases for $H$,   
where $C$ indicates the stability of the basis, see Section  2.1
for details.  These numbers are the main topic of our analysis.
We compare nonlinear approximations with linear approximations. 
Here we consider the class $\L_n$ of all continuous linear mappings 
$S_n : F \to H$, 
$$ 
S_n(f) = \sum_{i=1}^n L_i(f) \cdot \tilde  h_i 
$$ 
with arbitrary $\tilde h_i \in H$.  
The worst case error of optimal linear 
mappings is given by 
$$ 
e_n^\lin (S , F , H) = \inf_{S_n \in \L_n}  e(S_n, F,H) . 
$$ 
The third class of approximation methods  that we study in this paper
is the class of continuous mappings $\C_n$, given by arbitrary 
continuous mappings
$N_n : F \to \R^n$ and $\phi_n : \R^n \to H$. 
Again we define the worst case error of optimal continuous 
mappings by
$$ 
e_n^{\rm cont} (S, F , H) = \inf_{S_n \in \C_n}  e(S_n, F, H ) , 
$$ 
where $S_n = \phi_n \circ N_n$.  
These numbers, or slightly different numbers, were 
studied by different authors, cf. \cite{DHM89,DKLT93,DD96,Ma90}. 
Sometimes the $e_n^\cont$ are called 
manifold widths of $S$, see \cite{DKLT93}.

\begin{rem} {\rm 
\begin{itemize}
\item[i)]

A purpose of this paper is to compare the numbers 

$e_{n,C}^\non (S, F , H)$ with the numbers $e_n^\lin (S, F, H)$, 
where $S: F \to H$ is the restriction of $\A^{-1} : G \to H$ to 
$F \subset G$. 
In this sense we compare optimal linear approximation
of $S$ (i.e., by linear mappings of rank $n$)  with the 
best $n$-term approximation with respect to an optimal Riesz basis. 
\item[ii)] To avoid possible misunderstandings, it is important
to clarify the following point. In the realm of approximation theory, 
very often the term ``linear approximation''  is used for 
an approximation scheme that comes  
from a sequence of linear spaces that are 
uniformly refined, see, e.g., {\rm \cite{DV98}} for a detailed discussion.
However, in our definition of 
$e_n^\lin (S, F, H)$ we allow {\rm arbitrary} linear $S_n$, 
not only those that are based on uniformly refined subspaces. 
In this paper,  the latter will be denoted by  
{\rm uniform} approximation scheme.
\end{itemize} } 

\end{rem} 

For reader's convenience, we finish this section  by briefly summarizing the
main results of this paper. 

\begin{itemize} 
\item
Theorem 1:   Assume that $F \subset G$ is quasi-normed. 
Then 
$$ 
e_{n,C}^\non  (S, F , H) \, \ge \, \frac{1}{2C}  \, 
b_{m} (S, F , H)  
$$ 
holds for all $m \ge 4\, C^2 \, n $, where $b_n (S, F , H)$ 
denotes the  $n$-th Bernstein width of the operator $S$,
see Section 2.2 for details.
\item
Theorem 2 and Corollary 1: 
Assume that $F \subset G$ is a Hilbert space and 
$$ 
b_{2n} (S, F , H) \asymp b_n (S, F , H).
$$ 
Then 
$$ 
e_n^\lin (S, F , H) = e_n^\cont (S, F , H) 
\asymp e_{n,C}^\non  (S, F , H). 
$$ 
In this sense, approximation by linear mappings is as good as
approximation by nonlinear mappings. 
In this paper, `$a \asymp b$' always means
that both quantities can be uniformly bounded by a constant multiple of each
other. Likewise, `$\lsim$'  indicates inequality up to constant factors.
\item
Theorem 4: 
Assume that
$S: H^{-s}(\Omega)  \to H^{s}_0 (\Omega)$ is an isomorphism, 
with no further assumptions. 
Then we have for all $C \ge 1$
$$ 
e_n^\lin (S, H^{-s+t} , H^{s} ) \asymp 
e_{n,C}^\non (S, H^{-s+t} , H^{s} ) \asymp n^{-t/d} . 
$$
In this sense, approximation by linear mappings is as good as
approximation by nonlinear mappings.

\noindent 

Theorem 5: If we allow only function evaluations instead 
of general linear information, 
then the order of convergence drops down from 
$n^{-t/d}$ to $n^{(s-t)/d}$, where $t> s+ d/2$.
\item
In Theorem 6 and 7 we study the Poisson equation 
and the best $n$-term wavelet approximation. 
Theorem 6 shows that best $n$-term wavelet approximation might be 
suboptimal in general. 
Theorem 7, however, shows that for a polygonal domain in $\R^2$ 
best $n$-term wavelet approximation is almost optimal. 
\end{itemize}

Some of these results (Corollary 1, Theorem 4) might be surprising 
since there is a widespread believe that nonlinear approximation 
is better than approximation by linear operators. 
Therefore we  want to make the following remarks concerning our setting: 

\begin{itemize} 
\item
We allow arbitrary linear operators $S_n$ with rank $n$, not only 
those that are based on a uniform refinement. 
\item
We consider the worst case error with respect to the unit ball 
of a Hilbert space. 
\item 
Our results are concerned with  approximations, not with their 
numerical realization. 
For instance, the   construction   of an  optimal linear method
might require the  precomputation of a suitable basis 
(depending on $\A$),  which is usually
a prohibitive task. See also Remark 10,  where we discuss in more
detail why nonlinear approximation is very important for 
the approximation of elliptic problems. 
\item 
In another paper (in progress) we continue this work under the 
assumption that $F$ 
is a general Besov space. 
Then it turns out that for some parameters nonlinear approximation 
is essentially better than linear approximation. 
\end{itemize}

\section{Basic Concepts of Optimality} \label{basic} 

\subsection{Classes of Admissible Mappings}  

\subsubsection*{Nonlinear Mappings $S_n$} 

We will study certain approximations of $S$ based on Riesz bases, cf., e.g.,
Meyer \cite[page 21]{Me}.

\begin{definition}  \label{def1}   
Let $H$ be a Hilbert space. Then a sequence $h_1, h_2, \ldots $
of elements of $H$ is called a {\rm Riesz basis}  for $H$ if
there exist positive constants $A$ and $B$ such that, 
for every sequence of scalars
$\alpha_1, \alpha_2, \ldots \, $ with $\alpha_i \not= 0$ for 
only finitely many $i$,  we have
\begin{equation}\label{Riesz}
A \Big(\sum_{k} |\alpha_k|^2 \Big)^{1/2} \le \Big\|
\sum_{k} \alpha_k \, h_k  \Big\|_H \le B 
\Big( \sum_{k} |\alpha_k|^2 \Big)^{1/2} 
\end{equation}
and the vector space of finite sums $\sum \alpha_k \, h_k$
is dense in $H$.  
\end{definition}

\begin{rem} {\rm 
The constants $A,B$ reflect the stability of the basis. 
Orthonormal bases are those with
$A=B=1$. Typical examples of Riesz bases are the biorthogonal wavelet 
bases on $\Rd$ or on certain Lipschitz domains, 
cf. Cohen \cite[Sect. 2.6,~2.12]{C03}. } 
\end{rem}

In what follows 
\begin{equation}     \label{eq05} 
\B = \{ h_i \mid i \in \Nb \}
\end{equation} 
will always denote a Riesz basis of $H$ and $A$ and $B$ 
will be the corresponding 
optimal constants in (\ref{Riesz}). 
We study optimal approximations $S_n$ of $S=\A^{-1}$ of the form 
\begin{equation}     \label{eq06} 
S_n (f) = u_n = \sum_{k=1}^n c_k \, h_{i_k} , 
\end{equation} 
where $f=\A(u)$. 
Assuming that we can choose $\B$, we 
want to choose an optimal basis $\B$. 
What is the error of such an approximation $S_n$ and in which sense 
can we say that $\B$ and $S_n$ are optimal?

It is important to note that optimality of $S_n$ 
does not make sense for a single $u$: 
we simply can take a $\B$ where $h_1$ is a multiple 
of $u$,  and hence we can write the exact solution $u$ as 
$u_1=c_1 h_1$, i.e., with $n=1$. 
To define optimality of an approximation $S_n$ we need a 
suitable subset of $G$. We consider the worst 
case error 
\begin{equation}     \label{eq07} 
e(S_n, F,H) := 
\sup_{\Vert f \Vert_F \le 1} \Vert \A^{-1}(f)- S_n(f) \Vert_H , 
\end{equation} 
where $F$ is a normed (or quasi-normed) space, $F \subset G$. 
For a given basis $\B$ 
we consider the class $\N_n (\B)$ of all (linear or 
nonlinear) mappings of the form 
\begin{equation}     \label{eq10} 
S_n(f) = \sum_{k=1}^n c_k \,  h_{i_k} ,      
\end{equation} 
where the $c_k$ and the $i_k$ depend in an arbitrary way on $f$. 
Optimality is expressed by the quantity
\[
\sigma_n (\A^{-1}f,\B)_H :=
\inf_{i_1, \ldots , i_n}  \inf_{c_1, \ldots \, c_n} \| \A^{-1}(f) -  
\sum_{k=1}^n c_k \,  h_{i_k} \, \|_H \, .
\]
This reflects the best $n$-term approximation of $\A^{-1}(f)$. 
This subject is widely studied, see the surveys \cite{DV98}
and \cite{Tem03}.
Since $S_n$ is arbitrary, 
one  immediately obtains 
\begin{eqnarray*}
\inf_{S_n \in \N_n(\B)} 
\sup_{\| f\|_F \le 1} \| \A^{-1}(f) - S_n (f)\|_H &=&   
\sup_{\| f\|_F \le 1} 
 \inf_{S_n \in \N_n(\B)}  \| \A^{-1}(f) - S_n(f)\|_H \\
 &=&  
\sup_{\| f\|_F \le 1} \, \sigma_n (\A^{-1}f,\B)_H \, .  
\end{eqnarray*}

We allow the basis $\B$ to be chosen in a nearly arbitrary way. 
It is natural to assume some common stability 
of the bases under consideration.
For a real number $C \ge 1$ we define
\begin{equation}   \label{eq10a}  
\B_C := \Big\{\B : \, B/A \le C\Big\}.
\end{equation}  
We define the nonlinear widths
$e_{n,C}^\non (S,F,H)$ as 
\begin{equation}     \label{eq11} 
e_{n,C}^\non (S, F , H) = \inf_{\B \in \B_C} \inf_{S_n \in \N_n(\B)} 
e(S_n, F,H) . 
\end{equation} 
These numbers are the main topic of our analysis. 
They could be 
called the {\em errors of the best $n$-term approximation}
(with respect to the collection $\B_C$ of  Riesz basis of $H$), 
for brevity we call them {\em nonlinear widths}. 
In this paper we investigate the numbers $e^\non_{n,C} (S, F, H)$  
only when  $H$ is a Hilbert space.
More general concepts are introduced and investigated in \cite{Tem03}.

\begin{rem} {\rm  
It should be clear that the class $\N_n (\B)$ contains 
many mappings that are difficult to compute.
In particular, the number $n$ just reflects the dimension of a 
nonlinear manifold and has nothing to do with computational cost. 
Since we are interested 
in lower bounds, our results are strengthened by considering 
such a large class of approximations.  }  
\end{rem}

\begin{rem} {\rm  
It is obvious from the definition \eqref{eq11} that 
$S_n^* \in \N_n(\B)$ can be optimal for a given basis $\B$ 
in the sense that 
$$
e(S_n^* , F, H) \approx \inf_{S_n \in \N_n (\B)}  e(S_n, F, H) , 
$$
although the number
$e_{n,C}^\non (S, F , H)$ is much smaller, since the given $\B$ is 
far from being optimal. See also Remark~\ref{rem3}. }  
\end{rem}

\subsubsection*{Linear Mappings $S_n$}

Here we consider the class $\L_n$ of all continuous linear mappings 
$S_n : F \to H$, 
\begin{equation}     \label{eq08} 
S_n(f) = \sum_{i=1}^n L_i(f) \cdot \tilde  h_i 
\end{equation} 
with arbitrary $\tilde h_i \in H$.  
For each $S_n$ we define $e(S_n,F,H)$ by \eqref{eq07} and 
hence we can define the worst case error of optimal linear 
mappings by
\begin{equation}     \label{eq09} 
e_n^\lin (S , F , H) = \inf_{S_n \in \L_n}  e(S_n, F,H) . 
\end{equation} 
The numbers $e_n^\lin (S, F , H)$ (or slightly different numbers) 
are usually called 
{\it approximation numbers}  or  {\it linear widths}   of $S: F \to H$, 
cf. \cite{Ma90,Pi87,Pi85,Ti90}.

If $F$ is a space of functions on a set $\Omega$  such 
that function evaluation $f \mapsto f(x)$ is continuous, 
then one can define the {\it linear sampling numbers} 
\begin{equation}     \label{eq09a} 
g_n^\lin (S , F , H) = \inf_{S_n \in \L_n^{\rm std}}  e(S_n, F,H) , 
\end{equation} 
where $\L_n^{\rm std} \subset \L_n$ contains only 
those $S_n$ that are of the form
\begin{equation}     \label{eq08a} 
S_n(f) = \sum_{i=1}^n f(x_i) \cdot \tilde  h_i 
\end{equation} 
with $x_i \in \Omega$. 
For the numbers $g_n^\lin$ we only allow {\it standard information},
i.e., function values of the right-hand side. 
The inequality 
$g_n^\lin (S , F , H) \ge e_n^\lin (S, F, H)$ 
is trivial. 
One also might allow nonlinear $S_n = \phi_n \circ N_n$ with (linear) 
standard information 
$N_n(f) = (f(x_1), \dots , f(x_n))$ and arbitrary $\phi_n : \R^n \to
H$. This leads to the {\it sampling numbers} 
$g_n (S , F , H)$.

\subsubsection*{Continuous Mappings $S_n$}

Linear mappings $S_n$ are of the form $S_n = \phi_n \circ N_n$,  where 
both $N_n : F \to \R^n$ and $\phi_n : \R^n \to H$ are linear and 
continuous. If we drop the linearity condition,  then we obtain
the class of all continuous mappings $\C_n$, given by arbitrary 
continuous mappings
$N_n : F \to \R^n$ and $\phi_n : \R^n \to H$.
Again we define the worst case error of optimal continuous 
mappings by
\begin{equation}     \label{eq12} 
e_n^{\rm cont} (S, F , H) = \inf_{S_n \in \C_n}  e(S_n, F, H ) . 
\end{equation} 
These numbers, or variants of same, were 
studied by different authors, cf. \cite{DHM89,DKLT93,DD96,Ma90}. 
Sometimes these numbers are called 
manifold widths of $S$, see  \cite{DKLT93}. 
The inequalities 
\begin{equation}     \label{eq13} 
e_{n,C}^\non   (S, F , H) \le e_n^\lin (S, F , H)  
\end{equation} 
and
\begin{equation}  \label{eq13a} 
e_n^\cont  (S, F , H) \le e_n^\lin (S, F , H)
\end{equation} 
are,  of course,  trivial.

%&&&&&&&&&&&&&&&&&&&&&&&&&&&&&&&&&&&&&&&&&&&&&&&&&&&&&&&&&&&&&&&&&

%&&&&&&&&&&&&&&&&&&&&&&&&&&&&&&&&&&&&&&&&&&&&&&&&&&&&&&&&&&&&&&&&

\subsection{Relations to Bernstein Widths}   \label{sec2.4}  

%&&&&&&&&&&&&&&&&&&&&&&&&&&&&&&&&&&&&&&&&&&&&&&&&&&&&&&&&&&&&&&&&&

%&&&&&&&&&&&&&&&&&&&&&&&&&&&&&&&&&&&&&&&&&&&&&&&&&&&&&&&&&&&&&&&&&

The following quantities are useful for the understanding of 
$e_n^{\cont}$ and $e_n^{\non}$. 

\begin{definition}
The number 
$b_n (S, F , H)$, called the $n$-th {\rm Bernstein width}  of the operator
$S: \, F \to H$,
is the radius of the largest $(n+1)$-dimensional 
ball that is contained in $S(\{ \Vert f \Vert_F \le 1 \})$.
\end{definition}

\begin{rem}\label{Pietsch}{\rm 
In the literature there are used different definitions of Bernstein widths.
E.g. in Pietsch \cite{Pi74} the following version is given.
Let $X_n$ denote subspaces of $F$ of dimension $n$. Then
\[
\widetilde{b}_n (S,F,H) := \sup_{X_n \subset F} \inf_{x \in X_{n}, x\neq 0} \, 
\frac{\| Sx\|_H}{\| x\|_F}\,.
\] 
As long as $S$ is an injective mapping we obviously have
$b_{n} (S,F,H) = \widetilde{b}_{n+1}(S,F,H)$.}
\end{rem}

As it is well-known,  Bernstein widths are  useful for 
the proof of lower bounds, 
see \cite{DHM89,DD96,Pi85}.
The next lemma is certainly known. 
Since we could not find a reference, we include it with a proof.

\begin{lemma}\label{bernstein}
Let $n \in \Nb$ and assume that $F \subset G$ is quasi-normed. 
Then the inequality 
\begin{equation}     \label{eq14} 
b_n (S, F, H) \le  e_n^\cont (S, F, H)   %, d^n (S,F,H)\Big)
\end{equation} 
holds for all $n$. 
\end{lemma}

\begin{proof}
We assume that $S(\{ \Vert f \Vert_F \le 1 \} )$ contains an 
$(n+1)$-dimensional ball $B \subset H$ of radius $r$. 
We may assume that the center is in the origin. 
Let $N_n : F \to \R^n$ be continuous. 
Since $S^{-1}(B)$ is an $(n+1)$-dimensional bounded and symmetric 
neighborhood of 0, it follows from the Borsuk Antipodality Theorem, 
see \cite[par. 4]{D85},
that there exists an $f \in \partial S^{-1} (B)$ with 
$N_n (f)=N_n (-f)$ and hence 
$$
S_n(f) = \phi_n (N_n (f)) = \phi_n  (N_n (-f))= S_n(-f)
$$
for any mapping $\phi_n  : \R^n \to G$. 
Observe that $\Vert f \Vert_F=1$. 
Since  $\Vert S(f) - S(-f) \Vert =2r$ and 
$S_n(f)=S_n(-f)$,  we find that the maximal error of $S_n$ 
on $\{ \pm f \}$ is at least~$r$. 
This proves 
\[
b_n (S, F , H) \le e_n^\cont (S, F , H)\, . 
\]
\end{proof} 

We will see that the $b_n$ can also be used to prove lower bounds 
for the $e_{n,C}^\non$. 
As usual, $c_0$ denotes the Banach space of all sequences 
$x=(x_j)_{j=1}^\infty$
of real numbers such that $\lim_{j \to \infty} x_j =0$ and
equipped with the norm of $\ell_\infty$.

Lemma 2 below has a long history since it is central 
in the theory of $s$-numbers. 
See \cite[Lemma 2.9.6]{Pi87}, where also its use for proving
a result as Lemma 3 is exhibited. 

\begin{lemma} \label{l0}
Let $V$ denote an $n$-dimensional subspace of $ c_0$. 
Then there exists an element
$x \in V$ such that $\|x\|_\infty = 1$ and at 
least $n$ coordinates of $x= (x_1, x_2,\ldots)$
have absolute value $1$.
\end{lemma}

\begin{lemma}   \label{l1} 
Let $V_n$ be an $n$-dimensional subspace of the Hilbert space $H$.  
Let $\B$ be a Riesz basis with Riesz constants $0 < A \le B <\infty$. 
Then there is a nontrivial element  $x \in V_n$ such that 
$x = \sum_{j=1}^\infty x_j \, h_j $
and
\[
A \sqrt{n} \, \| (x_j)_j\|_\infty \le \| x \|_H . 
\] 
\end{lemma}

\begin{proof} 
Associated with any $x \in H$ there is a sequence $(x_j)_j$
of coefficients with respect to $\B$ that belongs to $c_0$.
In the same way,  we associate with $V_n \subset H$ a 
subspace $X_n \subset c_0$, also of dimension $n$.
As a consequence of Lemma \ref{l0},  we find 
an element $(x_j)_j \in X_n$ such that
\[
0 < |x_{j_1}| = \ldots = |x_{j_n}| = \| (x_j)_j \|_\infty  <\infty\, .
\] 
This implies
\[
\| x \|_H \ge A (\sum_{l =1}^n |x_{j_l}|^2)^{1/2} 
= A \sqrt{n} \, \| (x_j)_j \|_\infty \, . 
\]
\end{proof}

\begin{thm}   \label{t1} 
Assume that $F \subset G$ is quasi-normed. 
Then 
\begin{equation}     \label{eq15} 
e_{n,C}^\non  (S, F , H) \, \ge \, \frac{1}{2C}  \, 
b_{m} (S, F , H)  
\end{equation} 
holds for all $m \ge 4\, C^2 \, n $.
\end{thm}

\begin{proof} 
Let $\B$ be a Riesz basis with Riesz constants $A$ and $B$ and let $m>n$.
Assume that
$S(\{ \Vert f \Vert \le 1 \})$ contains an $m$-dimensional
ball with radius $\e$. Using  Lemma \ref{l1}, there exists an 
$x \in S(\{ \Vert f \Vert \le 1 \} )$ such that $x= \sum_i x_i \, h_i$,
$\Vert x \Vert = \e$ and $ |x_i| \le A^{-1} \, m^{-1/2} \e $ for all $i$.
Let $x_{1}, \ldots x_{n}$ be the $n$ largest components
(with respect to the absolute value) of $x$. 
Now, consider $y= \sum_i y_i \, h_i$ such that at most $n$
coefficients are nonvanishing. Then   
\[
\| x-y \|_H  \ge  A \| (x_i - y_i)_i\|_2 
 \]
and the optimal choice of $y$ (with respect to the right-hand side)
is given by $y^0$, where
$y_1^0 =x_1, \ldots \, , y_n^0 = x_n$.
Now we continue our estimate
\begin{eqnarray}\label{eq15b}
 A \| (x_i - y_i)_i\|_2  & \ge & A \, (\| (x_i)_i \|_2 - \| (y_i)_i \|_2) 
\nonumber
\\
& \ge & A \Big( \frac \e B - \frac 1A \, \e \, \sqrt{\frac{n}{m}}\Big) = 
\e \Big(\frac A B - \sqrt{\frac{n}{m}}\Big)\, .
\end{eqnarray}
The right-hand side is at least $ \e \, A/(2B)$ if $m \ge 4 B^2 n/A^2 $.
\end{proof}

\begin{rem} {\rm 
Probably the constant $1/(2C)$ is not optimal. But it is obvious from 
{\rm (\ref{eq15b})} that for $m$ tending to 
infinity the constant is approaching $A/B$. } 
\end{rem}

%&&&&&&&&&&&&&&&&&&&&&&&&&&&&&&&&&&&&&&&&&&&&&&&&&&&&&&&&&&&&&&&&&&

%&&&&&&&&&&&&&&&&&&&&&&&&&&&&&&&&&&&&&&&&&&&&&&&&&&&&&&&&&&&&&&&&&&

\subsection{The Case of a Hilbert Space}   \label{Hilbert}

%&&&&&&&&&&&&&&&&&&&&&&&&&&&&&&&&&&&&&&&&&&&&&&&&&&&&&&&&&&&&&&&&&&

%&&&&&&&&&&&&&&&&&&&&&&&&&&&&&&&&&&&&&&&&&&&&&&&&&&&&&&&&&&&&&&&&&&

%\subsection{General Facts on Hilbert Spaces}

%&&&&&&&&&&&&&&&&&&&&&&&&&&&&&&&&&&&&&&&&&&&&&&&&&&&&&&&&&&&&&&&&&&

%&&&&&&&&&&&&&&&&&&&&&&&&&&&&&&&&&&&&&&&&&&&&&&&&&&&&&&&&&&&&&&&&&&

Now let us assume, in addition to the assumptions 
of the previous subsections, that 
$F \subset G$ is a Hilbert space. 
The following result is well known, see \cite{Pi74}
and Remark~\ref{Pietsch}.

\begin{thm}   \label{t2} 
Assume that $F$ is a Hilbert space. Then 
\begin{equation}     \label{eq16} 
e_n^\lin (S, F , H) = e_n^\cont (S, F , H) = b_n (S, F , H). 
\end{equation} 
\end{thm} 

%  \begin{rem} {\rm  
%  There are different definitions of $s$-numbers in the literature. 
%  In particular, the definition used in the monograph \cite{Pi87} 
%  does not coincide with that one from the article \cite{Pi74}. 
%  In the sense of \cite{Pi87} the Bernstein numbers are not $s$-numbers,
%  but they are $s$-numbers in the sense of \cite{Pi74}.  }  
%  \end{rem} 
%  Arbitrary continuous mappings cannot be better than linear mappings. 
%  This is a general result for Hilbert spaces. 

In many applications one studies problems with ``finite smoothness'' 
and then, as a rule, one has the estimate 
\begin{equation}     \label{eq17} 
b_{2n} (S, F , H) \asymp b_n (S, F , H). 
\end{equation} 
Formula (\ref{eq17})  especially holds for the operator 
equations that we study in 
Section \ref{jetzt}. 
Then we conclude that approximation by 
optimal linear mappings yields the same order of convergence as
the best $n$-term approximation.

\begin{cor}     \label{cor1}  
Assume that $S: F \to H$ with Hilbert spaces $F$ and $H$,  with 
\eqref{eq17} holding. 
Then 
\begin{equation}     \label{eq18} 
e_n^\lin (S, F , H) = e_n^\cont (S, F , H) 
\asymp e_{n,C}^\non  (S, F , H). 
\end{equation} 
\end{cor}

%  \begin{proof} 
%  This follows from Theorem~\ref{t1} and Theorem~\ref{t2}, together 
%  with \eqref{eq17}.  
%  \end{proof} 

%&&&&&&&&&&&&&&&&&&&&&&&&&&&&&&&&&&&&&&&&&&&&&&&&&&&&&&&&&&&&&&&&&&

%&&&&&&&&&&&&&&&&&&&&&&&&&&&&&&&&&&&&&&&&&&&&&&&&&&&&&&&&&&&&&&&&&&

\section{Elliptic Problems} \label{jetzt}

%&&&&&&&&&&&&&&&&&&&&&&&&&&&&&&&&&&&&&&&&&&&&&&&&&&&&&&&&&&&&&&&&&&

%&&&&&&&&&&&&&&&&&&&&&&&&&&&&&&&&&&&&&&&&&&&&&&&&&&&&&&&&&&&&&&&&&&

In this section,  we study the more special case where 
$\Omega \subset \R^d$ is a bounded Lipschitz 
domain and
$\A=S^{-1}  : H^s_0 (\Omega) \to H^{-s} (\Omega)$  
is an isomorphism, where $s > 0$. 
The first step is  to recall the  definition of
the smoothness spaces that are needed for our analysis.

\subsection{Function Spaces}

If $m$ is a natural number,  we let  $H^m (\Omega)$
denote  the set of all 
functions $u \in L_2 (\Omega)$ such that the (distributional)
derivatives
$D^\alpha u$ of order $|\alpha| \le m$  also belong to $L_2 (\Omega)$.
This set,  equipped with the norm
\[
\| \, u \, \|_{H^m(\Omega)} := \sum_{|\alpha|\le m} 
\|\,  D^\alpha u\, \|_{L_2(\Omega)} , 
\]
becomes a Hilbert space. For a  positive noninteger $s$,  we 
define $H^s (\Omega)$ as specific Besov spaces. If $h\in\Rd$, we
let  $\Omega_h$ denote  the set of all $x\in\Omega$ 
such that the line segment
$[x,x+h]$ is contained in $\Omega$.  The modulus of smoothness
$\omega_r(u,t)_{L_p(\Omega)}$ of a 
function $u\in L_p(\Omega)$, where  $0<p\le \infty$,
is defined by  
$$ 
 \omega_r(u,t)_{L_p(\Omega)}:=\sup_{|h|\le
t}\|\Delta_h^r(u,\cdot)\|_{L_p(\Omega_{rh})},\quad t>0,
 $$
with $\Delta_h^r$ the $r$-th difference with step   $h$.
For $s>0$ and $0<q,p\le \infty$,  the {\em Besov space}  
$B_q^{s}(L_p(\Omega))$ is defined as the space of all 
functions $u \in L_p(\Omega)$ for which 
\begin{equation} \label{besovdef}
|u|_{B^{s}_q(L_p(\Omega))}:= \left\{ \begin{array}{ll}
\left(\int_0^{\infty}[t^{-s}\omega_r(u,t)_{L_p(\Omega)}]^qdt/t\right)^{1/q}, 
& 0<q<\infty,\\[1mm]
 \sup_{t\geq 0} t^{-s}\omega_r(u,t)_{L_p(\Omega)},    & 
q=\infty~,\end{array} \right. \end{equation}
is finite with $r\in \Nb$, $s < r \le s +1$, see, e.g., 
\cite{T92} for details.   
It turns out that   (\ref{besovdef}) is a (quasi-)semi-norm
for $B^{s}_q(L_p(\Omega))$. If we add $\|u\|_{L_p(\Omega)}$ to
(\ref{besovdef}), we obtain a (quasi-)norm for  
$B^{s}_q(L_p(\Omega))$.  Then, for positive
noninteger $s$,  we define
$$H^s(\Omega):=B^s_2(L_2(\Omega)).$$
It is known that this definition coincides up to 
equivalent norms with other definitions based,
e.g., on complex or real interpolation, cf. 
Dispa \cite{Di}, Lions and Magenes \cite[Vol.~1]{LM} and Triebel \cite{T02}.  

For all $s>0$ we let
$H^s_0 (\Omega)$ denote  the closure of the test
functions $\D (\Omega)$ in $H^s(\Omega)$. Finally, we put
\[
H^{-s} (\Omega) := (H^s_0 (\Omega))' \, , \qquad s >0\, \quad 
s \neq \frac 12 + k \, ,
\]
where $k \in \Nb_0$.
Alternatively (and this is done e.g. in \cite{JK1} and will play a role in
Subsection \ref{Poissonk})
one could use the following approach: define for $s>0$
\[%begin{eqnarray*}
\widetilde{H}^{s} (\Omega) :=   \Big\{ u \in L_2 (\Omega): ~~
\mbox{there exists} ~~ g \in   H^{s} (\R^d) ~~\mbox{with} ~~ 
g_{|_\Omega} = u ~~ \mbox{and}~~ \supp g \subset \overline{\Omega}\Big\}
\]%\end{eqnarray*}
equipped with the induced norm.
Then, for all $s>0$, $ s \neq \frac 12 + k$, $ k\in \Nb$  
\[
H^{s}_0 (\Omega) = \widetilde{H}^{s} (\Omega)\, , 
\]
in the sense of equivalent norms, cf. Grisvard \cite[Cor.~1.4.4.5]{Gr85}. 
If $0 < s = \frac 12 + k$, $ k\in \Nb$,  then we put
\begin{equation}\label{all}
H^{-s}(\Omega) := \Big(\widetilde{H}^{s} (\Omega)\Big)'\, .
\end{equation}
Observe, that by the previous remark this could be used as definition for all values of 
$s>0$ (up to equivalent norms).\\
Since we have Hilbert spaces, linear mappings are (almost) 
optimal approximations, i.e., Corollary~\ref{cor1} holds. 
We want to say more about the structure of an optimal linear $S_n$
for the approximation of $S=\A^{-1}$. For this, the  notion of a 
``regular problem'' is useful.

%&&&&&&&&&&&&&&&&&&&&&&&&&&&&&&&&&&&&&&&&&&&&&&&&&&&&&&&&&&&&&&&&&&

%&&&&&&&&&&&&&&&&&&&&&&&&&&&&&&&&&&&&&&&&&&&&&&&&&&&&&&&&&&&&&&&&&&

\subsection{Regular Problems}

%&&&&&&&&&&&&&&&&&&&&&&&&&&&&&&&&&&&&&&&&&&&&&&&&&&&&&&&&&&&&&&&&&&

%&&&&&&&&&&&&&&&&&&&&&&&&&&&&&&&&&&&&&&&&&&&&&&&&&&&&&&&&&&&&&&&&&&

The notion of regularity is very important for the theory and 
the numerical treatment of operator equations, see \cite{H92}. 
We use the following definition and assume that $t>0$.

\begin{definition} Let $s >  0$. 
An isomorphism $\A: H^s_0 (\Omega) \to H^{-s}(\Omega)$ 
is {\rm $H^{s+t}$-regular}
if also 
\begin{equation}     \label{eq19} 
\A: H^s_0(\Omega) \cap  H^{s+t} (\Omega) \to H^{-s+t} (\Omega)
\end{equation} 
is an isomorphism.  
\end{definition} 
A classical example is the Poisson equation
in a $C^{\infty}$-domain: this yields an operator that 
is $H^{1+t}$-regular for every $t>0$. 
We refer, e.g., to \cite{H92} for further information and examples. 
It is known that in this situation we obtain the optimal rate
\begin{equation} \label{rate} 
e_n^\lin (S, H^{-s+t}(\Omega)  , H^{s}(\Omega)  ) \asymp n^{-t/d} 
\end{equation} 
of linear methods. This is a classical result, at least 
for $t, s \in \Nb$ and for special domains.
We refer to the books \cite{Dy96,Pe96,We96} 
that contain hundreds of references.

We prove that the rate \eqref{rate} 
is true for arbitrary $s, t>0$, 
and for arbitrary bounded (nonempty, of course) Lipschitz domains. 
The optimal rate can be obtained by using 
Galerkin spaces that do not depend on the 
particular operator $\A$. 
With nonlinear approximations we cannot obtain a better rate 
of convergence.

\begin{thm}   \label{t3} 
Assume that the problem is $H^{s+t}$-regular.
Then for all $C \ge 1$, we have 
\begin{equation}     \label{eq20} 
e_n^\lin (S, H^{-s+t}(\Omega)  , H^{s}(\Omega)  ) \asymp 
e_{n,C}^\non (S, H^{-s+t}(\Omega)  , H^{s}(\Omega)  ) \asymp n^{-t/d} 
\end{equation} 
and the optimal order can be obtained by 
subspaces of $H^s$ that do not depend 
on the operator $S=\A^{-1}$. 
\end{thm}

\begin{proof} 
Consider first the identity (embedding) 
$I: H^{s+t}(\Omega)  \to H^s(\Omega) $. 
It is known that 
$$
e_n^\lin (I, H^{s+t}(\Omega) , H^s(\Omega) )  \asymp 
n^{-t/d}.
$$
This is a classical result (going back to Kolmogorov (1936), 
see \cite{K36}) for $s, t \in \Nb$, see also \cite{Pi85}.
For the general case ($s,t>0$ and arbitrary 
bounded Lipschitz domains)  see \cite{ET96} and \cite{T02}. 
We obtain the same order for 
$I: H^{s+t}(\Omega)\cap H^s_0 (\Omega)  \to H^s(\Omega) $.

We assume \eqref{eq19},  and hence 
$S: H^{-s+t} (\Omega) \to H^{s+t} (\Omega) \cap H^s_0 (\Omega)$ 
is an isomorphism. Hence we obtain the same order of the 
$e_n^\lin$ for $I$ and for $I \circ S_{|H^{-s+t} (\Omega)}$.
Together with  Corollary~\ref{cor1} this proves 
\eqref{eq20}.

Assume that the linear mapping
$$
\sum_{i=1}^n g_i \, L_i (f) 
$$
is good for the mapping $I: H^{s+t}(\Omega)  \cap H^s_0(\Omega)   
\to H^s(\Omega) $, 
i.e., we consider a sequence of such approximations
with the optimal rate $n^{-t/d}$. 
Then the linear mappings
$$
\sum_{i=1}^n g_i \, L_i (Sf) 
$$
achieve the optimal rate $n^{-t/d}$ for the mapping 
$S: H^{-s+t}(\Omega)  \to H^{s+t}(\Omega)  \subset H^s(\Omega) $. 
\end{proof}

\begin{rem}  {\rm  
The same $g_i$ are good for all $H^{s+t}(\Omega)$-regular 
problems on $H^{-s+t} (\Omega)$;
only the linear functionals, given by $L_i \circ S_{|H^{-s+k}}$, 
depend on the operator $\A$. 
For the numerical realization we can use the Galerkin method 
with the space $V_n$ generated by $g_1, \dots , g_n$. 
It is known that for $V_n$ one can take spaces that are based on
uniform refinement, e.g., constructed by
uniform grids or uniform finite elements schemes.
Indeed, if we consider
a sequence $V_n$ of uniformly refined spaces with
dimension $n$, then, under natural conditions, the following
characterization holds:
\begin{equation} \label{linchar} 
u \in H^{t+s}(\Omega) \Longleftrightarrow 
\sum_{n=1}^{\infty} [n^{t/d} E_n(u)]^2\frac{1}{n} < \infty, \quad \mbox
{where} \quad E_n(u):=\inf_{g\in V_n}\|u-g\|_{H^s}, \end{equation}
see, e.g,  {\rm \cite{DDD, Os}} and the references therein. } 
\end{rem}

\begin{rem}  {\rm  
Observe that the assumptions of Theorem~\ref{t3} are rather 
restrictive. Formally we assumed that $\Omega$ is an arbitrary bounded
Lipschitz domain and that $\A$ is $H^{s+t}$-regular. 
In practice, however, problems tend to be regular only if 
$\Omega$ has a smooth boundary. }  
\end{rem}

%&&&&&&&&&&&&&&&&&&&&&&&&&&&&&&&&&&&&&&&&&&&&&&&&&&&&&&&&&&&&&&&&&&

%&&&&&&&&&&&&&&&&&&&&&&&&&&&&&&&&&&&&&&&&&&&&&&&&&&&&&&&&&&&&&&&&&&

\subsection{Nonregular Problems}  \label{nonr}

%&&&&&&&&&&&&&&&&&&&&&&&&&&&&&&&&&&&&&&&&&&&&&&&&&&&&&&&&&&&&&&&&&&

%&&&&&&&&&&&&&&&&&&&&&&&&&&&&&&&&&&&&&&&&&&&&&&&&&&&&&&&&&&&&&&&&&&

The next result shows that linear approximations also give 
the optimal rate $n^{-t/d}$ in the nonregular case. 
An important difference,  however, is the fact that now 
the Galerkin space must depend on the operator $\A$. 
Related results can be found in the literature, see 
\cite{KS99,Me00,W87}. 
Again we allow arbitrary $s$ and $t>0$ and arbitrary bounded 
Lipschitz domains. We also prove that nonlinear 
approximation methods do not yield a better rate of convergence.

\begin{thm}   \label{t4} 
Assume that
$S: H^{-s}(\Omega)  \to H^{s}_0 (\Omega)$ is an isomorphism, 
with no further assumptions. 
Here $\Omega \subset \R^d$ is a bounded Lipschitz domain. 
Then we  have for all $C \ge 1$
\begin{equation}     \label{eq21} 
e_n^\lin (S, H^{-s+t}(\Omega)  , H^{s}(\Omega)  ) \asymp 
e_{n,C}^\non (S, H^{-s+t}(\Omega)  , H^{s}(\Omega)  ) \asymp n^{-t/d} . 
\end{equation} 
\end{thm}

\begin{proof} 
Consider first the identity (or embedding) 
$I: H^{-s+t}(\Omega)  \to H^{-s} (\Omega) $. 
It is known that 
$$
e_n^\lin (I, H^{-s+t}(\Omega) , H^{-s}(\Omega) )  \asymp 
n^{-t/d}.
$$
Again this is a classical result, for the general case 
(with $s, t>0$ and $\Omega$ an arbitrary 
bounded Lipschitz domain), see \cite{T02}. 

We assume that 
$S: H^{-s}(\Omega)  \to H^{s}_0 (\Omega)$ is an isomorphism, 
so that 
$e_n^\lin$ have the same order 
for $I$ and for $S \circ I$.
Together with Theorem~\ref{t1} and Corollary~\ref{cor1}, this proves 
\eqref{eq21}.

Assume that the linear mapping
$$
\sum_{i=1}^n g_i \, L_i (f) 
$$
is good for the mapping $I: H^{-s+t} \to H^{-s}$, 
i.e., we consider a sequence of such approximations 
with the optimal rate $n^{-t/d}$. 
Then the linear mappings
\begin{equation}     \label{eq22} 
\sum_{i=1}^n S(g_i)  \, L_i (f) 
\end{equation}  
achieve the optimal rate $n^{-t/d}$ for the mapping 
$S: H^{-s+t}(\Omega)  \to H^s (\Omega) $. 
\end{proof}

\begin{rem} {\rm    \label{remuni} 
It is well-known that uniform  methods can be quite bad for problems that
are not regular. Indeed, the  general characterization 
{\rm (\ref{linchar})} implies that
the approximation order  of uniform methods is 
determined by the Sobolev regularity of the solution $u$.
Therefore, if the problem is nonregular, i.e.,
if the solution $u$ lacks Sobolev smoothness, 
then the order of convergence of uniform methods
drops down.   }  
\end{rem}

\begin{rem}   {\rm    \label{rem3} 
For nonregular problems, we use linear combinations of $S(g_i)$. 
The $g_i$ do not depend on $S$, but of course the 
$S(g_i)$ do depend on $S$. 
This has important practical consequences: 
if we want to realize good approximations of the form 
\eqref{eq22} then we need to know the $S(g_i)$. 
Observe also that in this case we need good knowledge 
about the approximation of the embedding
$I : H^{-s+t}(\Omega)  \to H^{-s}(\Omega) $. For $s>0$,  
this embedding is not often studied in numerical analysis.

Hence we see an important difference between regular and 
arbitrary operator equations: Yes, the order of optimal 
linear approximations is the same in both cases and also 
nonlinear (best $n$-term) approximations cannot be better. 
But to construct good linear methods in the general case 
we have to know or to precompute the $S(g_i)$,   which is usually 
almost 
impossible in practice or at least much too expensive.

This leads us to the following problem: 
Can we find a $\B \in \B_C$ (here we think about a wavelet basis, 
but we do not want to exclude other cases) that depends 
only on $t$, $s$, and $\Omega$ such that 
\begin{equation}     \label{eq22a} 
\inf_{S_n \in \N_n (\B) } e(S_n, H^{-s+t} (\Omega) , H^s (\Omega) ) 
\asymp n^{-t/d} 
\end{equation} 
for many different operator equations,
given by an isomorphism $S=\A^{-1}: H^{-s} (\Omega) \to H_0^s (\Omega)$?

We certainly cannot expect that a single basis $\B$ is optimal for all 
reasonable operator equations, but the 
results in Section  \ref{Poissonk} indicate that wavelet
methods seem to have some potential in 
this direction. In any case it is important to 
distinguish between ``an approximation $S_n$ is optimal with respect 
to the given basis $\B$'' and ``$S_n$ is optimal with respect 
to the optimal basis $\B$''. 
See also {\rm \cite{Ni03}} and {\rm \cite{NW00}}. }  
\end{rem}

%%%%%%%%%%%%%%%%%%%%%%%%%%%%%%%%%%%%%%

\subsection{Function Values} \label{funct}

%%%%%%%%%%%%%%%%%%%%%%%%%%%%%%%%%%%%%%%%

Now we study the numbers 
$g_n(S, H^{t-s} (\Omega)  , H^s (\Omega) ) = g_n^\lin (S, H^{t-s}
(\Omega) , H^s (\Omega) )$ 
under similar conditions as we had in Theorem~\ref{t4}. 
In particular we do not assume that the problem is regular. 
However we have to assume $t > s + d/2$ so that 
function values will 
continuously depend on $f \in H^{t-s} (\Omega) $.

Consider first the embedding
$I: H^{t}(\Omega)  \to L_2 (\Omega) $, where 
$\Omega$ is a bounded Lipschitz domain $\Omega \subset \R^d$. 
We want to use function values of $f \in H^t (\Omega)$ 
and hence have to assume that $t > d/2$.  
It is known that 
\begin{equation}  \label{fu1} 
e_n^\lin (I, H^{t}(\Omega) , L_2 (\Omega) )  \asymp 
g_n^\lin (I, H^t (\Omega) , L_2 (\Omega)) \asymp 
n^{-t/d},
\end{equation} 
see \cite{NT04}. 
This means that arbitrary linear functionals do not yield 
a better order of convergence than function values. 
Observe that we always have 
$g_n=g_n^\lin$, since we consider mappings between Hilbert 
spaces and hence the linear spline algorithm is always 
optimal, see \cite[4.5.7]{TWW88}. 

It is interesting that for $s>0$ arbitrary linear information 
is superior to function evaluation.
In the theorem that follows, we make 
no smoothness or regularity assumptions. 

\begin{thm}   \label{t4a} 
Assume that
$S: H^{-s}(\Omega)  \to H^{s}_0 (\Omega)$ is an isomorphism, 
where $\Omega \subset \R^d$ is a bounded Lipschitz domain. 
Then
\begin{equation}     \label{eq21a} 
g_n (S, H^{-s+t} (\Omega)  , H^{s} (\Omega)  ) = 
g_n^\lin  (S, H^{-s+t} (\Omega)  , H^{s} (\Omega)  )
\asymp n^{(s-t)/d} , 
\end{equation} 
for $t > s+ d/2$.  
\end{thm}

\begin{proof} 
As in the proof of 
Theorem~\ref{t4}, it is enough to prove that 
\begin{equation}    \label{eq21b}  
g_n (I, H^{-s+t}(\Omega) , H^{-s}  (\Omega) ) \asymp
n^{(s-t)/d}.
\end{equation}  
To prove the upper and the lower bound for \eqref{eq21b}, 
we use several auxiliary problems and 
start with the upper bound. It is known from \cite{NT04} that 
$$
g_n (I, H^{-s+t}(\Omega) , L_2 (\Omega) )  \asymp n^{(s-t)/d}.
$$
provided that $t-s>d/2$.
{}From this we obtain the upper bound
$$
g_n (I, H^{-s+t}(\Omega) , H^{-s}  (\Omega) ) \le c \cdot 
n^{(s-t)/d}
$$
by embedding.

For the lower bound we use the bound 
\begin{equation}   \label{eq21c} 
g_n (I, H^{-s+t}(\Omega) , L_1 (\Omega) )  \asymp n^{(s-t)/d}, 
\end{equation}  
again from \cite{NT04}. The  lower bound in \eqref{eq21c} 
is proved by the technique of bump functions: 
Given $x_1 , \dots , x_n \in \Omega$, one can construct a function 
$f \in H^{-s+t} (\Omega) $ with norm one such that 
$f(x_1) = \dots = f(x_n) =0$ and  
\begin{equation}   \label{eq21d}  
\Vert f \Vert_{L_1} \ge c \cdot n^{(s-t)/d} ,
\end{equation}  
where $c>0$ does not depend on the $x_i$ or on $n$. 
The same technique can be used to prove lower bounds for integration
problems. 
We consider an integration problem 
\begin{equation}  \label{eq21e}  
{\rm Int} (f) = \int_\Omega f \sigma \, dx ,
\end{equation}  
where $\sigma   \ge 0$ is a smooth (and nonzero) function 
on $\Omega$ with compact support. 
Then this technique gives: 
Given $x_1 , \dots , x_n \in \Omega$, one can construct a function 
$f \in H^{-s+t} (\Omega) $ with norm one such that 
$f(x_1) = \dots = f(x_n) =0$ and  
\begin{equation}   \label{eq21f}  
{\rm Int} (f) \ge c \cdot n^{(s-t)/d} ,
\end{equation}  
where $c>0$ does not depend on the $x_i$ or on $n$. 
Since we assumed that $\sigma$ is smooth with compact support, 
we have
$$
\Vert f \Vert_{H^{-s}} \ge c \cdot | {\rm Int} (f) | 
$$
and hence we may replace in \eqref{eq21f} ${\rm Int} f$ by $\Vert f 
\Vert_{H^{-s}}$, hence 
$$
g_n (I, H^{-s+t} (\Omega), H^{-s} (\Omega) ) \ge c \cdot n^{(s-t)/d} .
$$
\end{proof}

%&&&&&&&&&&&&&&&&&&&&&&&&&&&&&&&&&&&&&&&&&&&&&&&&&&&&&&&&&&&&&&&&&&

%&&&&&&&&&&&&&&&&&&&&&&&&&&&&&&&&&&&&&&&&&&&&&&&&&&&&&&&&&&&&&&&&&&

\subsection{The Poisson Equation} \label{Poissonk}

%&&&&&&&&&&&&&&&&&&&&&&&&&&&&&&&&&&&&&&&&&&&&&&&&&&&&&&&&&&&&&&&&&&

%&&&&&&&&&&&&&&&&&&&&&&&&&&&&&&&&&&&&&&&&&&&&&&&&&&&&&&&&&&&&&&&&&&

Finally we discuss our results for the specific case of
the Poisson equation 
\begin{eqnarray}
-\triangle u &=&f \quad  \mbox{in}\quad \Omega \label{Poisson}\\
           u&=&0 \quad \mbox{on} \quad \partial \Omega \nonumber
\end{eqnarray}
on a bounded Lipschitz domain $\Omega$ contained in $\R^d$, $d\ge 2$. 
Here, as always in this paper, we understand Lipschitz domain in the sense of
Steins notion of domains with minimal smooth boundary, cf. Stein 
\cite[VI.3]{St70}.

It is well-known that (\ref{Poisson}) fits 
into our setting with $s=1$. Indeed,
if we consider the weak formulation of this problem, it can be checked that
(\ref{Poisson}) induces a boundedly invertible 
operator $\A=\triangle: ~H_0^1(\Omega)\longrightarrow H^{-1}(\Omega)$,
see again \cite[Chapter 7.2]{H92} for details. 

%E  Naechster Satz weg. 

%  Here we meet the problem of existence 
%  of an appropriate  Riesz basis for $H^1_0 (\Omega)$. 
In this section, we shall
especially focus on {\em wavelet} 
bases $\Psi=\{\psi_{\lambda} : \lambda \in {\cal J}\}.$ 
The indices $\lambda \in  {\cal J}$ typically encode several
types of information, namely the {\em scale}  (often denoted $|\lambda|$),
the spatial location and also the type of the wavelet. 
Recall that in a classical
setting,  a tensor product construction 
yields $2^d-1$ types of wavelets \cite{Me}.  For instance,
on the real line $\lambda$ can be identified with $(j,k)$, 
where $j=|\lambda|$ denotes the dyadic
refinement level and $2^{-j}k$ signifies the location of the  wavelet.
We will not discuss at this point any technical description of the basis
$\Psi$. Instead we assume that the domain $\Omega$ under consideration
enables us to construct a wavelet basis $\Psi$ with the following properties:
\begin{itemize}
\item the wavelets are {\em local} in the sense that
$$ 
\mbox{diam}(\mbox{supp} \psi_{\lambda}) \asymp  2^{-|\lambda|}, 
\quad \lambda \in {\cal J};
$$
\item the wavelets satisfy the {\em cancellation property}
$$
|\langle v, \psi_{\lambda}\rangle| \lsim 2^{-|\lambda|{\widetilde{m}}}
\|v\|_{H^{\widetilde{m}}({\rm supp} \, \psi_{\lambda})},$$
where $\widetilde{m}$ denotes some suitable parameter, and 
\item 
the wavelet basis induces characterizations of Besov spaces of the form
\begin{equation} \label{besovchar} 
\|f\|_{B^s_q(L_p(\Omega))} 
\asymp \left(\sum_{|\lambda|=j_0}^{\infty}
2^{j(s+d(\frac{1}{2}-\frac{1}{p}))q}\left(\sum_{\lambda \in {\cal J}, 
|\lambda|=j}|\langle f,
\tilde{\psi}_{\lambda}\rangle|^p\right)^{q/p}\right)^{1/q},
\end{equation}
where $s>d\left(\frac{1}{p} -1\right)_{+}$ 
and
$\tilde{\Psi}=\{\tilde{\psi}_{\lambda} 
: \lambda \in {\cal J}\}$ denotes the
 {\em dual basis}
$$
\langle \psi_{\lambda}, 
\tilde{\psi}_{\nu}\rangle =\delta_{\lambda, \nu}, \quad \lambda, \nu
\in {\cal J}.
$$
\end{itemize}

%  \fix{E: Fuer welche $p,q$ usw. brauchen wir das? 
%  (Frage von Triebel)}
%  \fix{S: Bemerkung eingefuegt}

For the applications we have in mind, especially the case 
$p=q, p \leq 2,  1/p \leq s/d +1/2$ is important, see, e.g., 
Theorem \ref{opti1} for details. 

By  exploiting  the norm equivalence (\ref{besovchar}) 
and using the fact that
$B^s_2(L_2(\Omega))=H^s(\Omega)$, a simple rescaling 
immediately yields  a Riesz basis for $H^s$.  
We shall also assume that the Dirichlet
boundary conditions can be included,  
so that a  characterization of the form(\ref{besovchar})
also  holds  for 
$H^s_0(\Omega)$. We refer to \cite{C03} for a 
% \fix{Seitenzahl oder so?}
detailed  
discussion.  In this setting, the following theorem holds. 

\begin{thm} \label{opti1} 
Let $S$ denote the solution operator for the problem {\rm (\ref{Poisson})}.  
Then, for sufficiently large $C$, 
best $n$-term wavelet approximation $S_n$ yields
% \fix{E: $S_n$ eingefuegt} 
\begin{eqnarray*} %\label{toll}
e^\non_{n,C} (S, H^{t-1}(\Omega) , H^1(\Omega) ) & \le & 
e(S_n, H^{t-1}(\Omega) , H^1(\Omega) ) \\
& \leq & c \, \left\{\begin{array}{lll}
n^{-\frac{t}{d}+\varepsilon} &\qquad & \mbox{if}\quad 0 <t\le 1/2\, , \\
&& \\
n^{-\frac{(t+1)}{3d}+ \varepsilon}&\qquad & \mbox{if}\quad \frac 1 2 < t
\le  \frac{d+2}{2(d-1)} \, , \\
&& \\
n^{-\frac{d}{2d(d-1)}+\varepsilon}
&\qquad & \mbox{if}\quad  \frac{d+2}{2(d-1)} < t \, , \\
\end{array}
\right.
\end{eqnarray*}
where $\varepsilon >0$ is arbitrary and $c$ does not depend on $n\in \Nb$. 
\end{thm}

\begin{proof} 
{\em Step 1.} All what we need from the wavelet basis
is the following estimate for the best $n$-term approximation
in the $H^1$-norm:
\begin{equation} \label{H1scale}
\|\, u-S_n(f)\, \|_{H^1} \leq c\,  
|u|_{B^{\alpha}_{\tau*}(L_{\tau*}(\Omega))}\, 
n^{(\alpha-1)/d}, \qquad 
\frac{1}{\tau^*}= \frac{(\alpha-1)}{d} + \frac{1}{2},
\end{equation}
see, e.g., \cite{DDD} for details.  
We therefore have to estimate
the Besov norm $B^{\alpha}_{\tau*}(L_{\tau*}(\Omega))$. 
\\
{\em Step 2.} Besov regularity of $u$.\\
First of all, we estimate  the Besov norm of $u$ in the specific
scale 
\begin{equation} \label{nonscale}
B^s_{\tau}(L_{\tau}(\Omega)), \quad
\hbox{where} \  \  \frac{1}{\tau} =\frac{s}{d}+\frac{1}{2}.
\end{equation}
Regularity estimates  in the scale (\ref{nonscale}) have already
been performed in \cite{DD}. We write the solution $u$ to (\ref{Poisson})
as $$u= \tilde{u} + v,$$
where $\tilde{u}$ solves $-\triangle \tilde{u} = \tilde{f}$ on a
smooth domain $\wt{\Omega} \supset \Omega$. Here  $\tilde{f}= {\mathcal E}(f)$
where ${\mathcal E}$ denotes some suitable extension operator
with respect to $\Omega$.  Furthermore, $v$
is the solution to the additional homogeneous Dirichlet problem
\begin{eqnarray} \triangle v&=& 0 \qquad \mbox{in}\quad\Omega
\label{dirichlet} \\
 v&=& g= - \mbox{Tr}(\tilde{u})\qquad \mbox{on}\quad \partial \Omega\, 
.\nonumber
\end{eqnarray}
%  \fix{S:  a priori Abschaetzung wieder eingefuegt}
{\em Substep 2.1} Regularity of $\tilde{u}$.
%  \fix{S:  Referenz geaendert}
Let $t>0$.
Let ${\mathcal E}$ be a bounded linear extension operator from
$B^{t-1}_2 (L_2 (\Omega)) \to B^{t-1}_2 (L_2 (\R^d))$, 
cf. \cite{Ry99}.
Then, by classical elliptic regularity on smooth domains,  
cf. e.g. \cite[Thm.~0.3]{JK1},
it follows  from ${\mathcal E} f \in  B^{t-1}_2 (L_2 (\widetilde{\Omega})) $ that
 $\tilde{u} \in  B^{t+1}_2 (L_2 (\widetilde{\Omega})) $ and 
$$ 
\|\, \tilde{u}\, \|_{B^{t+1}_2(L_2(\wt{\Omega}))} 
\leq c_1 \, \|{\mathcal {E}}\|\, \,   \| f \|_{B^{t -1}_2(L_2(\Omega))}\, .
$$
Known embeddings of Besov spaces yield 
\begin{equation}\label{a2}
\|\, \tilde{u}\, \|_{B^{t +1-\varepsilon}_{q}(L_{q}(\wt{\Omega}))}\leq 
c_2 \, \|\, \tilde{u}\, \|_{B^{t +1}_2(L_q(\wt{\Omega}))} \leq
c_3 \, \|\, \tilde{u}\, \|_{B^{t +1}_2(L_2(\wt{\Omega}))} \leq c_4 \| f \|_{B^{t -1}_2(L_2(\Omega))}\, ,  
\end{equation}
where $0 < q \le 2$ and $\varepsilon >0$ are arbitrary.
\\
{\em Substep 2.2} The regularity of $v$.
An important theorem of Jerison and Kenig, see  \cite{JK2,JK3} and also
\cite[Thm.~5.1]{JK1} (for $d\ge 3$), reads as 
\begin{equation}\label{a1}
\|\, v\, \|_{B_2^{\rho}(L_2(\Omega))} 
\leq  c_5 \, \|\, g\, \|_{B^{\rho-1/2}_2(L_2(\partial \Omega))}, 
\quad \mbox{if} \quad 1/2 < \rho  < 3/2\, .
\end{equation}
Trace problems for Lipschitz boundaries are investigated in \cite{JW}.
We refer to this monograph and to \cite{JK1} also for the exact meaning 
of $\mbox{Tr}$ and $B^\rho_2 (L_2 (\partial \Omega))$, respectively. 
Theorem 2 on page 209 in \cite{JW} and (\ref{a1}) yield
\begin{eqnarray*} 
\|\, v \, \|_{B^{\rho}_2(L_2(\Omega))} \le c_5\, 
\|\,  \mbox{Tr}\, \| 
\, \| \, \tilde{u}\, \|_{B^{\rho}_2(L_2(\wt \Omega))}
\leq c_6 \,  \|\, \tilde{u}\, \|_{B^{t +1}_2(L_2(\wt{\Omega}))} \leq c_7 \| f \|_{B^{t -1}_2(L_2(\Omega))}\,, 
\end{eqnarray*}
if $1/2 < \rho < 3/2$ and $\rho \le t+1$. Consequently
$v \in B^\vartheta_2 (L_2 (\Omega))$ with $\vartheta < \min (3/2, t+1)$.
A harmonic function in a bounded Lipschitz domain has a higher Besov
regularity than Sobolev regularity. More precisely,
\[
v \in B^s_\tau (L_\tau (\Omega))\, , \qquad 0 < s < 
\frac{\vartheta d}{d-1}\, , \quad \frac{1}{\tau} = \frac sd + \frac 12\, ,  
\]
and 
\[
\| \, v \, \|_{B^s_\tau (L_\tau (\Omega))} \le c_8 \, 
\| \, v\, \|_{B^\vartheta_2 (L_2 (\Omega))}\, , 
\]
see \cite{DD}.
Combining this with (\ref{a2}) we arrive at
\begin{equation}\label{a3}
u \in B^s_\tau (L_\tau (\Omega))\, , \qquad 0 < s < \min \Big( 
\frac{\vartheta d}{d-1}, t+1 \Big)\, , 
\quad \frac{1}{\tau} = \frac sd + \frac 12\, ,  
\end{equation}
together with the estimate
$$
\|u\|_{B^{s}_{\tau}(L_{\tau}(\Omega))}\leq
\|\tilde{u}\|_{B^{s}_{\tau}(L_{\tau}(\Omega))} 
+
\|v\|_{B^{s}_{\tau}(L_{\tau}(\Omega))}\\
\leq  c_9\|f\|_{B^{t-1}_2(L_2(\Omega))}.
$$
{\em Substep 2.3} An interpolation argument.
Another theorem of Jerison and Kenig, see \cite[Thm.~0.5]{JK1},
yields 
\begin{equation}\label{a4}
u \in B^s_2 (L_2 (\Omega))\, , \quad s< \min (3/2, t+1) \,, \quad \mbox{where}\quad \|u\|_{B^s_2(L_2(\Omega))}\leq c_{10}\|f\|_{B^{t-1}_2(L_2(\Omega))}. 
\end{equation}
Thanks to the real interpolation formula
\[
\Big( B^{s_0}_{p_0} (L_{p_0} (\Omega)),  B^{s_1}_{p_1} (L_{p_1} (\Omega)) 
\Big)_{\Theta,p} = B^s_p (L_p (\Omega))\, , 
\qquad \mbox{(equivalent quasi-norms)}, 
\]
\[
0 < \Theta < 1 \, , \quad 
s= (1-\Theta)\,  s_0 + \Theta \, s_1\, , \quad
\frac 1p = \frac{1-\Theta}{p_0} + \frac{\Theta}{p_1}\, , 
\]
valid for all $s_0,s_1 \in {\mathbb R}$ and all $0 < p_0,p_1 < \infty$,
%\fix{W: $p < \infty$ wegen Referenz}
cf. \cite{T02}, we can combine these two different assertions 
(\ref{a3}), (\ref{a4}) about the regularity of $u$.
Let e.g. $1/2 \le t \le (d+2)/(2d-2) $.
Then we use the interpolation formula with
\[
s_0 = 3/2-\varepsilon\, , \qquad s_1 = t+1\, ,
\qquad \mbox{and}\qquad \Theta = 1/3\, , 
\]  
and find that $u \in B^s_{\tau^*} (L_{\tau^*} (\Omega))$, 
$s= 1+ (t+1)/3- \varepsilon'$, $1/\tau^* = (s-1)/d + 1/2 $,
where 
$$
\|u\|_{B^{s}_{\tau^*}(L_{\tau*}(\Omega))}
\leq c_{11} \|f\|_{B^{t -1}_{2}(L_2(\Omega))},$$
and  $\varepsilon'$ can be made as small as we want.  In summary, we have
$$ 
\sup_{\|f\|_{B^{t-1}_2(L_2(\Omega))}\leq 1}\|u-S_n(f)\|_{H^1} 
\leq c_{12}\,  n^{-(\frac{t+1}{3d}+\varepsilon)}.
$$
The other two cases can be treated in an analoguous way. We omit details. 
\end{proof}

Theorem \ref{opti1} shows that  best $n$-term wavelet approximation 
might be suboptimal in 
general. However, for more specific domains, i.e., 
for polygonal domains, much more can be said.
Let $\Omega$  denote  a simply connected polygonal 
domain contained in $\R^2$, the segments of $\partial \Omega$ are
denoted
by $\overline{\Gamma}_l, \Gamma_l$ open, $l=1, \ldots, N$ 
numbered in positive  orientation. Furthermore, $\Upsilon_l$
denotes the endpoint of $\Gamma_l$ and $\omega_l$ denotes 
the measure of the interior  angle at $\Upsilon_l$. 
Then the following theorem holds:

\begin{thm} \label{optipol}
Let $S$ denote the solution operator for the problem {\rm (\ref{Poisson})}
in a polygonal domain in $\R^2$. Let $k$ be a nonnegative integer such that
\[
k \neq \frac{m\pi}{\omega_l} \qquad \mbox{for all} \quad  m \in \Nb, \quad l=1, \ldots\, , N\, .
\]
Then, for sufficiently large $C$, 
best $n$-term wavelet approximation $S_n$ yields 
\begin{equation}\label{toller}
e^\non_{n,C} (S, H^{k-1}(\Omega) , H^{1}(\Omega) )  \le  
e(S_n, H^{k-1}(\Omega) , H^{1}(\Omega) )  \leq C \, n^{-k/2 + \varepsilon}\, , 
\end{equation}
where $\epsilon>0$ is arbitrary and $C$ does not depend on $n$. 
 \end{thm}
%\fix{W: gaendert} 

 \begin{proof}  
The proof is based on the fact that $u$ can be decomposed 
into a regular part $u_R$ and a singular part $u_S$,
 $u=u_R + u_S,$ where $u_R \in B^{k +1}_2(L_2(\Omega))$ 
and $u_S$  only depends on the shape of 
the domain and can be computed explicitly. 
This result was established by Grisvard, 
see \cite{Gr2} or \cite[Chapt.~4,~5]{Gr85}, and  \cite[Sect.~2.7]{Gr1} for details. 
%  \fix{W: Referenz ist jetzt so richtig, war ein Tippfehler von mir} 
We introduce
polar coordinates $(r_l, \theta_l)$ in the 
vicinity of each vertex $\Upsilon_l$ and introduce the functions
$$ 
{\cal S}_{l,m}(r_l,\theta_l) = 
\zeta_l(r_l)r_l^{\lambda_{l,m}}\sin(m\pi \theta_l /\omega_l), 
$$
when
$\lambda_{l,m}:=m\pi/\omega_l$
is not an integer and
$$
{\cal S}_{l,m}(r_l,\theta_l) = 
\zeta_l(r_l)r_l^{\lambda_{l,m}}[\log r_l \sin(m\pi \theta_l/\omega_l)+
\theta_l\cos(m\pi \theta_l/\omega_l)]
$$
otherwise, 
$m \in {\mathbb N}$, $l=1, \ldots \, N$.  
Here $\zeta_l$ denotes a suitable $C^{\infty}$ truncation function. 
Then for $f \in H^{k-1} (\Omega)$ one has
  \begin{equation} \label{singpart}
  u_S= \sum_{l=1}^N \sum_{0<\lambda_{l,m}< k} \, c_{l,m} \, {\cal S}_{l,m},
  \end{equation}
  provided that  no $\lambda_{l,m}$ is equal to $k$. This means that
the finite number of singularity functions that is 
needed depends on the scale of  spaces
  we are interested in, i.e., on the smoothness parameter $k$. 
  According to  (\ref{H1scale}), we have to estimate the
  Besov regularity of both, $u_S$ and $u_R$, in the specific scale
  $$
  B^{\alpha}_{\tau^*}(L_{\tau^*}(\Omega))  
  \qquad \frac{1}{\tau^*}= \frac{(\alpha-1)}{d} + \frac{1}{2}.  $$
Since $u_R \in B_2^{k+1}(L_2(\Omega))$,  
classical embeddings of Besov spaces imply that
\begin{equation} \label{regreg}
u_R \in {B^{k+1-\epsilon'}_{\tau*}(L_{\tau*}(\Omega))}
\qquad \frac{1}{\tau^*}= \frac{(k-\epsilon')}{d} 
+ \frac{1}{2}\qquad \mbox{for arbitray small } \quad \epsilon' >0.\end{equation}
Moreover, it has been shown in \cite{Da1} that the 
functions  ${\cal S}_{l,m}$ defined above satisfy 
\begin{equation} \label{singreg} 
{\cal S}_{l,m}(r_l,\theta_l)\in {B^{\alpha}_{\tau*}(L_{\tau*}(\Omega))},  
\qquad \frac{1}{\tau^*}= \frac{(\alpha-1)}{d} +
\frac{1}{2}\quad 
\mbox{for all} \quad \alpha >0.\end{equation}
By combining (\ref{regreg}) and (\ref{singreg}) we see that 
\[
u \in 
{B^{k+1-\epsilon'}_{\tau*}(L_{\tau*}(\Omega))}
\qquad \frac{1}{\tau^*}= \frac{(k-\epsilon')}{d} 
+ \frac{1}{2}\qquad \mbox{for arbitray small } \quad \epsilon' >0.
\]
To derive an estimate uniformly with respect to the unit ball in $H^{k-1} (\Omega)$
we argue as follows. We put
%  \fix{W: l-Abhängigkeit angegeben}
\[
{\mathcal N} := \span \Big\{ {\cal S}_{l,m}(r_l,\theta_l): \quad 0 < \lambda_{m,l}< k\, , 
\: l=1,\,  \ldots \, , N\Big\}\, .
\]
Let $\gamma_l$ be the trace operator with respect to the segment $\Gamma_l$.
Grisvard has shown that
$\Delta$ maps 
\[
H := \Big\{ u \in H^{k+1}(\Omega): \quad \gamma_l u =0\, , \, l = 1, \ldots \, ,N\Big\} \: + \: {\mathcal N} 
\]
onto $H^{k-1}(\Omega)$, cf. \cite[Thm.~5.1.3.5]{Gr85}. This mapping is also injective, see
 \cite[Lemma~4.4.3.1, Rem.~5.1.3.6]{Gr85}.
We equip the space $H$ with the norm
%  \fix{W: j zu l gemacht}
\[
\| \, u\,  \|_H := \| \, u_R + u_S\,  \|_H =
\| \, u_R\,  \|_{H^{k+1}(\Omega)} + \sum_{l=1}^N \sum_{0<\lambda_{l,m}< k} \, |c_{l,m}|\, , 
\] 
see (\ref{singpart}). Then it becomes a Banach space. Furthermore, $\Delta$ is continuous.
Banach's continuous inverse theorem implies  that the solution operator is continuous considered 
as a mapping  from $H^{k-1}(\Omega)$ onto $H$.
Observe
%  \fix{W: j zu l gemacht}
\[
\| \, u_R + u_S \, \|_{B^{k+1-\epsilon'}_{\tau*}(L_{\tau*}(\Omega))}  \le
C\, \Big(
\| \, u_R \, \|_{B^{k+1}_{2}(L_{2}(\Omega))} + \sum_{l=1}^N \sum_{0<\lambda_{l,m}< k} \, |c_{l,m}|\Big) 
\]
with some constant $C$ independent of $u$. 
\end{proof}

\medskip
\noindent
{\bf Acknowledgment. } \
We thank Stefan Heinrich, Aicke Hinrichs, 
Hans Triebel, Art Werschulz and two referees for their 
valuable remarks and comments.

%&&&&&&&&&&&&&&&&&&&&&&&&&&&&&&&&&&&&&&&&&&&&&&&&&&&&&&&&&&&&&&&&&&

%&&&&&&&&&&&&&&&&&&&&&&&&&&&&&&&&&&&&&&&&&&&&&&&&&&&&&&&&&&&&&&&&&&&&

\bigskip
\vbox{\noindent Stephan Dahlke\\
Philipps-Universit\"at Marburg\\
FB12 Mathematik und Informatik\\
Hans-Meerwein Stra\ss e\\
Lahnberge\\
35032 Marburg\\
Germany\\
e--mail: {\tt dahlke@mathematik.uni-marburg.de}\\
WWW: {\tt http://www.mathematik.uni-marburg.de/$\sim$dahlke/}\\}

\bigskip
\smallskip
\vbox{\noindent Erich Novak, Winfried Sickel\\
Friedrich-Schiller-Universit\"at Jena\\ 
Mathematisches Institut\\
Ernst-Abbe-Platz 2\\ 
07743 Jena \\ 
Germany\\
e-mail: {\tt \{novak, sickel\}@math.uni-jena.de}\\
WWW: {\tt http://www.minet.uni-jena.de/$\sim$\{novak,sickel\}/}\\}

\bigskip

\end{document}